\documentclass[12pt,leqno]{article}

\def\tr{\mathop{\rm tr}}

\newtheorem{theorem}{Theorem}
\newtheorem{lemma}[theorem]{Lemma}
\newtheorem{proposition}[theorem]{Proposition}
\newtheorem{definition}[theorem]{Definition}
\newtheorem{remark}[theorem]{Remark}

\newcommand{\begintheorem}{\addtocounter{equation}{1}\begin{theorem}}
\newcommand{\beginlemma}{\addtocounter{equation}{1}\begin{lemma}}
\newcommand{\beginproposition}{\addtocounter{equation}{1}\begin{proposition}}
\newcommand{\begindefinition}{\addtocounter{equation}{1}\begin{definition}}
\newcommand{\beginremark}{\addtocounter{equation}{1}\begin{remark}}

\begin{document}

\title{Potpourri, 2}

\author{Stephen William Semmes	\\
	Rice University		\\
	Houston, Texas}

\date{}

\maketitle


\renewcommand{\thefootnote}{}   

\footnotetext{These notes are connected to the ``potpourri''
topics class in the Department of Mathematics, Rice University,
in the fall semester of 2004.}

\tableofcontents

\bigskip

	As usual, ${\bf R}$ and ${\bf C}$ denote the real and complex
numbers.  If $z$ is a complex number, then $z$ can be expressed as $x
+ y \, i$, where $x$, $y$ are real numbers, called the real and
imaginary parts of $z$, and $i^2 = -1$.  In this case the complex
conjugate of $z$ is denoted $\overline{z}$ and defined to be $x - y \,
i$.

	If $x$ is a real number, then the absolute value of $x$ is
denoted $|x|$ and defined to be equal to $x$ when $x \ge 0$ and to
$-x$ when $x \le 0$.  Notice that $|x + y| \le |x| + |y|$ and $|x \,
y| = |x| \, |y|$ for all real numbers $x$, $y$.

	If $z = x + y \, i$ is a complex number, $x, y \in {\bf R}$,
then the modulus of $z$ is denoted $|z|$ and is the nonnegative real
number defined by $|z|^2 = z \, \overline{z} = x^2 + y^2$.  For every
pair of complex numbers $z$, $w$ we have that $|z + w| \le |z| + |w|$
and $|z \, w| = |z| \, |w|$.

\section{Normed vector spaces}
\label{normed vector spaces}
\setcounter{equation}{0}

	Let $V$ be a real or complex vector space.  By a \emph{norm}
on $V$ we mean a function $\|v\|$ defined for all $v \in V$ such that
$\|v\|$ is a nonnegative real number for all $v \in V$ which is equal
to $0$ if and only if $v = 0$,
\begin{equation}
	\|\alpha \, v\| = |\alpha| \, \|v\|
\end{equation}
for all real or complex numbers $\alpha$, as appropriate, and all $v
\in V$, and
\begin{equation}
	\|v + w\| \le \|v\| + \|w\|
\end{equation}
for all $v, w \in V$.  In this event $\|v - w\|$ defines a metric on
$V$.

	An \emph{inner product} on a real or complex vector space $V$
is a function $\langle v, w \rangle$ defined for $v, w \in V$ which
takes values in the real or complex numbers, according to whether $V$
is a real or complex vector space, and which satisfies the following
properties.  First, for each fixed $w \in V$, $\langle v, w \rangle$
is linear as a function of $v$, which is to say that
\begin{equation}
	\langle v + v', w \rangle 
		= \langle v, w \rangle + \langle v', w \rangle
\end{equation}
for all $v, v' \in V$ and
\begin{equation}
	\langle \alpha \, v, w \rangle = \alpha \, \langle v, w \rangle
\end{equation}
for all real or complex numbers $\alpha$, as appropriate, and all $v
\in V$.  Second, $\langle v, w \rangle$ is symmetric in $v$, $w$ when
$V$ is a real vector space, which means that
\begin{equation}
	\langle w, v \rangle = \langle v, w \rangle
\end{equation}
for all $v, w \in V$, and it is Hermitian-symmetric in the complex
case, which means that
\begin{equation}
	\langle w, v \rangle = \overline{\langle v, w \rangle}
\end{equation}
for all $v, w \in V$.  As a result, $\langle v, w \rangle$ is linear
in $w$ for each fixed $v$ in the real case, and it is conjugate-linear
in the complex case.  From the symmetry condition it follows that
$\langle v, v \rangle$ is a real number for all $v \in V$ even in the
complex case, and the third condition is that
\begin{equation}
	\langle v, v \rangle \ge 0
\end{equation}
for all $v \in V$, with
\begin{equation}
	\langle v, v \rangle > 0
\end{equation}
when $v \ne 0$.

	Suppose that $V$ is a real or complex vector space with inner
product $\langle v, w \rangle$, and put
\begin{equation}
	\|v\| = \langle v, v \rangle^{1/2}
\end{equation}
for $v \in V$.  The Cauchy--Schwarz inequality states that
\begin{equation}
	|\langle v, w \rangle| \le \|v\| \, \|w\|
\end{equation}
for all $v, w \in V$.  This can be derived using the fact that
\begin{equation}
	\langle v + \, \alpha w, v + \, \alpha w \rangle
\end{equation}
is a nonnegative real number for all scalars $\alpha$.  As a
consequence, one can check that
\begin{equation}
	\|v + w\| \le \|v\| + \|w\|
\end{equation}
for all $v, w \in V$, by expanding $\|v + w\|^2$ in terms of the inner
product.  It follows that $\|v\|$ defines a norm on $V$.

	A subset $E$ of a vector space $V$ is said to be \emph{convex}
if for each $v, w \in E$ and each real number $t$ with $0 \le t \le 1$
we have that
\begin{equation}
	t \, v + (1-t) \, w \in E.
\end{equation}
Suppose that $\|v\|$ is a real-valued function on $V$ such that $\|v\|
\ge 0$ for all $v \in V$, $\|v\| = 0$ if and only if $v = 0$, and
$\|\alpha \, v\| = |\alpha| \, \|v\|$ for all scalars $\alpha$ and all
$v \in V$.  Then $\|v\|$ defines a norm on $V$ if and only if
\begin{equation}
	\{v \in V : \|v\| \le 1 \}
\end{equation}
is a convex subset of $V$.  In other words, this is equivalent to the
triangle inequality in the presence of the other conditions.  This is
not too difficult to verify.

	Fix a positive integer $n$, and consider ${\bf R}^n$, ${\bf
C}^n$ as real or complex vector spaces.  More precisely, ${\bf R}^n$
and ${\bf C}^n$ consist of $n$-tuples of real or complex numbers, as
appropriate.  If $v = (v_1, \ldots, v_n)$, $w = (w_1, \ldots, w_n)$
are elements of ${\bf R}^n$ or of ${\bf C}^n$, then the sum $v + w$ is
defined coordinatewise,
\begin{equation}
	v + w = (v_1 + w_1, \ldots, v_n + w_n).
\end{equation}
Similarly, if $\alpha$ is a real or complex number and $v = (v_1,
\ldots, v_n)$ is an element of ${\bf R}^n$ or of ${\bf C}^n$, as
appropriate, then the scalar product $\alpha \, v$ is defined
coordinatewise,
\begin{equation}
	\alpha \, v = (\alpha_1 \, v_1, \ldots, \alpha_n \, v_n).
\end{equation}

	The standard inner product on ${\bf R}^n$, ${\bf C}^n$ is
defined by
\begin{equation}
	\langle v, w \rangle = \sum_{j=1}^n v_j \, w_j
\end{equation}
in the real case and
\begin{equation}
	\langle v, w \rangle \sum_{j=1}^n v_j \, \overline{w_j}
\end{equation}
in the complex case.  If $p$ is a real number with $1 \le p < \infty$,
then we put
\begin{equation}
	\|v\|_p = \biggl(\sum_{j=1}^n |v_j|^p \biggr)^{1/p}
\end{equation}
in both the real and complex cases, and we can extend this to $p =
\infty$ by
\begin{equation}
	\|v\|_\infty = \max (|v_1|, \ldots, |v_n|).
\end{equation}
Notice that $\|v\|_2$ is the norm associated to the standard inner
product on ${\bf R}^n$, ${\bf C}^n$.  For all $p$, $1 \le p \le
\infty$, one can check that $\|v\|_p$ defines a norm on ${\bf R}^n$,
${\bf C}^n$.  The triangle inequality is easy to check when $p = 1,
\infty$, and when $1 < p < \infty$ one can show that the closed unit
ball associated to $\|v\|_p$ is convex using the fact that $t^p$ is a
convex function on the nonnegative real numbers.

	Suppose that $1 \le p, q \le \infty$ and that
\begin{equation}
	\frac{1}{p} + \frac{1}{q} = 1,
\end{equation}
with $1/\infty = 0$, in which case we say that $p$, $q$ are
``conjugate exponents''.  If $v$, $w$ are elements of ${\bf R}^n$ or
of ${\bf C}^n$, then
\begin{equation}
	\biggl| \sum_{j=1}^n v_j \, w_j \biggr| \le \|v\|p \, \|w\|_q.
\end{equation}
This is H\"older's inequality.

	When $p = q = 2$ H\"older's inequality reduces to the
Cauchy--Schwarz inequality.  When $p, q = 1, \infty$ one can check it
directly using the triangle inequality for scalars.  Now suppose that
$1 < p, q < \infty$, and observe that
\begin{equation}
	a \, b \le \frac{a^p}{p} + \frac{b^q}{q}
\end{equation}
for all nonnegative real numbers $a$, $b$, as a result of the
convexity of the exponential function, for instance.  Hence
\begin{equation}
	\biggl| \sum_{j=1}^n v_j \, w_j \biggr|
		\le \frac{\|v\|_p^p}{p} + \frac{\|w\|_q^q}{q}
\end{equation}
for all $v$, $w$ in ${\bf R}^n$ or in ${\bf C}^n$, by applying the
previous inequality to $|v_j \, w_j|$ and summing over $j$.  This
yields H\"older's inequality when $\|v\|_p = 1$ and $\|w\|_q = 1$, and
the general case follows from a scaling argument.

	The triangle inequality for $\|v\|_p$ is known as Minkoski's
inequality, and one can also derive it from H\"older's inequality, in
analogy with the $p = 2$ case.  Let us restrict our attention to $1 <
p < \infty$, since the $p = 1, \infty$ cases can be handled directly.
For all $v$, $w$ in ${\bf R}^n$ or in ${\bf C}^n$ we have that
\begin{equation}
	\|v + w\|_p^p \le \sum_{j=1}^n |v_j| \, |v_j + w_j|^{p-1}
		+ \sum_{j=1}^n |w_j| \, |v_j + w_j|^{p-1}.
\end{equation}
If $q$ is the exponent conjugate to $p$, then H\"older's inequality
implies that
\begin{equation}
	\|v + w\|_p^p 
  \le (\|v\|_p + \|w\|_p) 
	\, \biggl(\sum_{j=1}^n |v_j + w_j|^{q(p-1)} \biggr)^{1/q}.
\end{equation}
This can be rewritten as
\begin{equation}
	\|v + w\|_p^p \le (\|v\|_p + \|w\|_p) \, \|v + w\|_p^{p-1},
\end{equation}
which implies that $\|v + w\|_p \le \|v\|_p + \|w\|_p$, as desired.

	If $v$ is an element of ${\bf R}^n$ or of ${\bf C}^n$ and $1
\le p < \infty$, then
\begin{equation}
	\|v\|_\infty \le \|v\|_p.
\end{equation}
More generally, if $1 \le p \le q < \infty$, then
\begin{equation}
	\|v\|_q \le \|v\|_p.
\end{equation}
Indeed,
\begin{equation}
	\|v\|_q^q = \sum_{j=1}^n |v_j|^q 
		\le \|v\|_\infty^{q-p} \, \|v\|_p^p \le \|v\|_p^q.
\end{equation}
Of course
\begin{equation}
	\|v\|_p \le n^{1/p} \, \|v\|_\infty
\end{equation}
for all $v$ in ${\bf R}^n$ or ${\bf C}^n$ and $1 \le p < \infty$.  For
$1 \le p \le q < \infty$ one can check that
\begin{equation}
	\|v\|_p \le n^{(1/p) - (1/q)} \, \|v\|_q
\end{equation}
using H\"older's inequality.

	If $V$ is a real or complex vector space and $\|v\|$ is a norm
on $V$, then
\begin{equation}
	\|v\| \le \|w\| + \|v - w\|
\end{equation}
and
\begin{equation}
	\|w\| \le \|v\| + \|v - w\|
\end{equation}
for all $v, w \in V$, by the triangle inequality.  Therefore
\begin{equation}
	\Bigl| \|v\| - \|w\| \Bigr| \le \|v - w\|
\end{equation}
for all $v, w \in V$.  In particular, $\|v\|$ is continuous on $V$ as
a real-valued function with respect to the metric associated to the
norm on $V$.  

	If $V$ is ${\bf R}^n$ or ${\bf C}^n$, then it is easy to see
that $\|v\|$ is bounded by a constant times the standard Euclidean
norm $\|v\|_2$, by expressing $v$ as a linear combination of the
standard basis vectors.  It follows that $\|v\|$ is continuous as a
real-valued function with respect to the usual Euclidean topology on
${\bf R}^n$ or ${\bf C}^n$.  By standard results from advanced
calculus, the minimum of $\|v\|$ over the set of $v$'s such that
$\|v\|_2 = 1$ is attained, since the latter is compact, and of course
the minimum is positive because $\|v\| > 0$ when $v \ne 0$.  This
implies that $\|v\|$ is also greater than or equal to a positive
constant times $\|v\|_2$.  As a consequence, the topology determined
by the metric $\|v - w\|$ is the same as the standard Euclidean
topology on ${\bf R}^n$, ${\bf C}^n$, as appropriate.

\section{Separation theorems}
\label{separation theorems}
\setcounter{equation}{0}

	Fix a positive integer $n$, and let $E$ be a nonempty closed
convex subset of ${\bf R}^n$.  Also let $p$ be a point in ${\bf R}^n$
which is not in $E$.  There exists a point $q \in E$ such that the
Euclidean distance $\|p - q\|_2$ from $p$ to $q$ is as small as
possible.

	Let $H$ be the affine hyperplane through $q$ which is
orthogonal to $p - q$.  In other words, using the standard inner
product on ${\bf R}^n$, $H$ consists of the $v \in {\bf R}^n$ such
that the inner product of $v - q$ with $p - q$ is equal to $0$.

	If $x$ is any element of $E$, then $x$ lies in the closed
half-space in ${\bf R}^n$ which is bounded by $H$ and which does not
contain $p$.  This is equivalent to saying that the inner product of
$x - q$ with $p - q$ is less than or equal to $0$, while the inner
product of $p - q$ with itself is equal to $\|p - q\|_2^2 > 0$.  One
can see this through a simple geometric argument, to the effect that
if $x \in E$ lies in the open half-space in ${\bf R}^n$ containing
$p$, then there is a point along the line segment joining $x$ to $q$
which is closer to $p$ than $q$ is.

	It follows that in fact $E$ is equal to the intersection of
the closed half-spaces containing it.  Namely, each point in ${\bf
R}^n$ which is not in $E$ is also in the complement of one of the
closed half-spaces containing $E$.

\beginremark {\rm The use of the Euclidean norm here may seem a bit
strange, since the statement that $q \in E$ and $H$ is a hyperplane
through $q$ such that $E$ is contained in a closed half-space in ${\bf
R}^n$ bounded by $H$ and $p$ is in the complementary open half-space
bounded by $H$ does not require the Euclidean norm or inner product.
One could just as well use a different inner product on ${\bf R}^n$,
which could lead to a different choice of $q$ and $H$.  Observe
however that if $q$ and $H$ the properties just mentioned, then there
is an inner product on ${\bf R}^n$ such that the distance from $q$ to
$p$ in the corresponding norm is as small as possible and $H$ is the
hyperplane through $q$ which is orthogonal to $q - p$.}
\end{remark}

	Next suppose that $E$ is a closed convex subset of ${\bf R}^n$
and that $p$ is a point in the boundary of $E$.  Thus $p \in E$ and
there is a sequence of points $\{p_j\}_{j=1}^\infty$ in ${\bf R}^n
\backslash E$ which converges to $p$.  In this case there is a
hyperplane $H$ in ${\bf R}^n$ which passes through $p$ such that $E$
is contained in one of the closed half-spaces bounded by $H$.  We can
reformulate this by saying that there is a vector $v \in {\bf R}^n$
such that $\|v\|_2 = 1$ and for each $x \in E$ we have that the inner
product of $x - p$ with $v$ is greater than or equal to $0$.  From the
previous argument we know that for each $j$ there is a point $q_j \in
E$ such that $\|q_j - p_j\|_2$ is as small as possible and for each $x
\in E$ the inner product of $x - q_j$ with $q_j - p$ is greater than
or equal to $0$.

	Put $v_j = (q_j - p) / \|q_j - p\|_2$, so that $\|v_j\|_2 = 1$
for all $j$ by construction.  By passing to a subsequence if necessary
we may assume that $\{v_j\}_{j=1}^\infty$ converges to a vector $v \in
{\bf R}^n$ such that $\|v\|_2 = 1$.  It is easy to check that $v$ has
the required properties, since $\{q_j\}_{j=1}^\infty$ converges to
$p$.

	Now let $\mathcal{C}$ be a closed convex cone in ${\bf R}^n$,
which means that $\mathcal{C}$ is a closed subset of ${\bf R}^n$, $0
\in \mathcal{C}$, for each $v \in \mathcal{C}$ and positive real
number $t$ we have that $t \, v \in \mathcal{C}$, and for each $v, w
\in \mathcal{C}$ we have that $v + w \in \mathcal{C}$.  Suppose that
$z$ is a point in ${\bf R}^n$ which is not in $\mathcal{C}$, so that
$t \, z$ is not in $\mathcal{C}$ for any positive real number $t$.
Let us check that there is a hyperplane in ${\bf R}^n$ which passes
through $0$ such that $\mathcal{C}$ is contained in one of the closed
half-spaces in ${\bf R}^n$ bounded by $H$ and $z$ is contained in the
complementary open half-space bounded by $H$.  This is equivalent to
saying that there is a vector $v \in {\bf R}^n$ such that $\|v\|_2 =
1$, the inner product of $v$ with any element of $\mathcal{C}$ is
greater than or equal to $0$, and the inner product of $v$ with $z$ is
strictly less than $0$.

	From the earlier arguments there is a $q \in \mathcal{C}$ such
that for any $x \in \mathcal{C}$ the inner product of $x - q$ with $q
- z$ is greater than or equal to $0$.  Put $v = (q - z) / \|q -
z\|_2$, so that $\|v\|_2 = 1$ automatically and the inner product of
$x - q$ with $v$ is greater than or equal to $0$ for all $x \in
\mathcal{C}$.  Because $0$ and $2 \, q$ are elements of $\mathcal{C}$,
we have that the inner product of $-q$, $q$ with $v$ are greater than
or equal to $0$, which is to say that the inner product of $q$ with
$v$ is actually equal to $0$.  Thus the inner product of any $x \in
\mathcal{C}$ with $v$ is greater than or equal to $0$, and the inner
product of $z$ with $v$ is negative is equal to the inner product of
$z - q$ with $v$, which is $- \|z - q\|_2 < 0$.

\section{Dual spaces}
\label{dual spaces}
\setcounter{equation}{0}

	Let $V$ be a finite-dimensional real or complex vector space,
and let $V^*$ denote the corresponding dual vector space of linear
functionals on $V$.  Thus the elements of $V^*$ are linear mappings
from $V$ into the real or complex numbers, as appropriate.  If $V$ is
a real or complex vector space, then $V^*$ is too, because one can add
linear functionals and multiply them by scalars.

	Suppose that the dimension of $V$ is equal to $n$, and that
$v_1, \ldots, v_n$ is a basis for $V$.  Thus every element of $V$ can
be expressed in a unique manner as a linear combination of the
$v_j$'s.  If $\lambda \in V^*$, then $\lambda$ is uniquely determined
by the $n$ scalars $\lambda(v_1), \ldots, \lambda(v_n)$, and these
scalars may be chosen freely.  Thus $V^*$ also has dimension equal to
$n$.

	Suppose that $V$ is equipped with a norm $\|v\|$.  If
$\lambda$ is a linear functional on $V$, then there is a nonnegative
real number $k$ such that
\begin{equation}
	|\lambda(v)| \le k \, \|v\|
\end{equation}
for all $v \in V$.  To see this, one might as well assume that $V$ is
equal to ${\bf R}^n$ or ${\bf C}^n$, using a basis for $V$ to get an
isomorphism with ${\bf R}^n$ or ${\bf C}^n$.  As we have seen, any
norm on ${\bf R}^n$ or ${\bf C}^n$ is equivalent to the standard
Euclidean norm, in the sense that each is bounded by a constant
multiple of the other.  The existence of an $k \ge 0$ as above then
follows from the corresponding statement for the Euclidean norm.

	Let us define the dual norm of a linear functional $\lambda$
on $V$ associated to the norm $\|v\|$ on $V$ by
\begin{equation}
	\|\lambda\|_* = \sup \{|\lambda(v)| : v \in V, \|v\| \le \}.
\end{equation}
Equivalently,
\begin{equation}
	|\lambda(v)| \le \|\lambda\|_* \, \|v\|
\end{equation}
for all $v \in V$, and $\|\lambda\|_*$ is the smallest nonnegative
real number with this property.  One can check that $\|\lambda \|_*$
does indeed define a norm on $V^*$.

	For instance, suppose that $V$ is equipped with an inner
product $\langle v, w \rangle$.  For each $w \in V$,
\begin{equation}
	\lambda_w(v) = \langle v, w \rangle
\end{equation}
defines a linear functional on $V$, and in fact every linear
functional on $V$ arises in this manner.  With respect to the norm on
$V$ associated to the inner product, the dual norm of $\lambda_w$ is
less than or equal to the norm of $w$, by the Cauchy--Schwartz
inequality.  By choosing $v = w$ one can check that the dual norm of
$\lambda_w$ is equal to the norm of $w$.

	Now suppose that $V$ is ${\bf R}^n$ or ${\bf C}^n$, and for
each $w$ in ${\bf R}^n$ or ${\bf C}^n$, as appropriate, consider the
linear functional $\lambda_w$ on $V$ given by
\begin{equation}
	\lambda_w(v) = \sum_{j=1}^n v_j \, w_j,
\end{equation}
$v = (v_1, \ldots, v_n)$, $w = (w_1, \ldots, w_n)$.  Every linear
functional on $V$ arises in this manner.  If $V$ is equipped with the
norm $\|v\|_p$ from Section \ref{normed vector spaces}, $1 \le p \le
\infty$, then the dual norm of $\lambda_w$ is equal to $\|w\|_q$,
where $q$ is the exponent conjugate to $p$.  Indeed, the dual norm of
$\lambda_w$ is less than or equal to $\|w\|_q$ by H\"older's
inequality.  Conversely, one can show that the dual norm of
$\lambda_w$ is greater than or equal to $\|w\|_q$ through specific
choices of $v$.

	If $V$ is any finite-dimensional real or complex vector space
equipped with a norm $\|\cdot \|$, and if $v$ is any vector in $V$,
then
\begin{equation}
	|\lambda(v)| \le \|\lambda\|_* \, \|v\|,
\end{equation}
just by the definition of the dual norm.  It turns out that for each
nonzero vector $v \in V$ there is a linear functional $\lambda$ on $V$
such that $\|\lambda\|_* = 1$ and $\lambda(v) = \|v\|$.  To see this
we may as well assume that $\|v\| = 1$, by scaling.

	Assume first that $V$ is a real vector space, which we may as
well take to be ${\bf R}^n$.  The closed unit ball in ${\bf R}^n$
associated to $\|\cdot \|$, consisting of vectors with norm less than
or equal to $1$, is a compact convex subset of ${\bf R}^n$, and $v$
lies in the boundary of this convex set, since $\|v\| = 1$.  As in
Section \ref{separation theorems}, there is a hyperplane $H$ through
$v$ such that the closed unit ball associated to the norm is contained
in one of the closed half-spaces bounded by $H$.  The linear
functional $\lambda$ on ${\bf R}^n$ that we want is
characterized by
\begin{equation}
	H = \{x \in {\bf R}^n : \lambda(x) = 1\}.
\end{equation}

	Now suppose that $V$ is a complex vector space.  If $\lambda$
is a linear mapping from $V$ to the complex numbers, then the real
part of $\lambda$ is a real-linear mapping from $V$ into the real
numbers, i.e., it is a linear functional on $V$ as a real vector
space, without the additional structure of scalar multiplication by
$i$.  Conversely, if one starts with a real-linear mapping from $V$
into the real numbers, then that is the real part of a unique
complex-linear mapping from $V$ into the complex numbers, as one can
verify.

	If $\beta$ is any complex number, then the modulus of $\beta$
can be described as the supremum of the real part of $\alpha \,
\beta$, where $\alpha$ runs through all complex numbers with $|\alpha|
\le 1$.  As a result, if $V$ is a complex vector space, $\|\cdot \|$
is a norm on $V$, and $\lambda$ is a complex linear functional on $V$,
then the dual norm of $\lambda$ can be described equivalently as the
supremum of the real part of $\alpha \, \lambda (z)$ as $\alpha$ runs
through all complex numbers with $|\alpha| \le 1$ and $z$ runs through
all vectors in $V$ with $\|z\| \le 1$.  Hence the dual norm of
$\lambda$ is equal to the supremum of the real part of $\lambda(z)$ as
$z$ runs through all vectors in $V$ with $\|z\| \le 1$, because one
can absorb the scalar factors into $z$.  In other words, the norm of
$\lambda$ as a complex linear functional on $V$ is equal to the norm
of the real part of $\lambda$ as a real linear functional on $V$,
using the same norm on $V$.  This permits one to derive the complex
case of the statement under consideration from the real case.

\section{Quotient spaces, norms}
\label{quotient spaces, norms}
\setcounter{equation}{0}

	Let $V$ be a finite-dimensional real or complex vector space,
and let $W$ be a linear subspace of $V$.  Consider the quotient space
$V / W$, which is basically defined by identifying points in $V$ whose
difference lies in $W$.  There is a canonical quotient mapping $q$,
which is a linear mapping from $V$ onto $W$.

	Suppose also that $V$ is equipped with a norm $\|v\|$.  Thus
we get a metric associated to this norm, and with respect to this
metric $W$ is a closed subset of $V$.  Indeed, in ${\bf R}^n$ or ${\bf
C}^n$, every linear subspace is a closed subset.  The general case can
be derived from this one because $V$ is isomorphic to ${\bf R}^n$ or
${\bf C}^n$ for some $n$, and the norm on $V$ is equivalent to the
usual Euclidean norm on ${\bf R}^n$ or ${\bf C}^n$, as appropriate,
with respect to this isomorphism.

	Let us define a quotient norm $\|\cdot \|_{V / W}$ on $V / W$
by saying that the norm of a point in $V / W$ is equal to the infimum
of the norms of the points in $V$ which are identified to that point
in the quotient.  Equivalently, for each $x \in V$, the norm of $q(x)$
in $V / W$ is equal to the infimum of the norms of $x + w$ in $V$ over
$w \in W$.  In particular, $\|q(x)\|_{V / W} \le \|x\|$ for all $x \in
V$.  It is not difficult to check that this does indeed define a norm
on the quotient space $V / W$.

	Now let $Z$ be a linear subspace of $V$, and suppose that
$\lambda$ is a linear functional on $Z$.  We would like to extend
$\lambda$ to a linear functional on $V$ whose dual norm on $V$ is the
same as that of $\lambda$ on $Z$, using the restriction of the given
norm on $V$ to $Z$ as a norm on $Z$.  We may as well assume that
$\lambda$ is not the zero linear functional, which is to say that
$\lambda(z) \ne 0$ for at least some $z \in Z$.  Let $W$ denote the
kernel of $\lambda$, which is the linear subspace of $Z$ consisting of
all vectors $y \in Z$ such that $\lambda(y) = 0$.

	Let us work in the quotient space $V / W$.  The quotient $Z /
W$ is a one-dimensional subspace of $V / W$.  Because $W$ is the
kernel of $\lambda$, there is a canonical linear functional on $Z / W$
induced by $\lambda$.  One can check that the norm of this linear
functional on $Z / W$, associated to the quotient norm on $Z / W$
obtained from our original norm $\|\cdot \|$ on $V$, is equal to the
norm of $\lambda$ as a linear functional on $Z$.

	Because $Z / W$ has dimension equal to $1$, there is a linear
functional $\mu$ on $V / W$ with dual norm equal to $1$ with respect
to the quotient norm on $V / W$ such that for each element of $Z / W$,
the absolute value or modulus of $\mu$ applied to that element is
equal to the quotient norm of that element.  This follows from the
result discussed in Section \ref{dual spaces}, applied to $V / W$ with
the dual norm.  We can multiply $\mu$ by a scalar to get a linear
functional on $V / W$ whose norm is equal to the norm of $\lambda$ on
$Z$ and which agrees on $Z / W$ with the linear functional induced
there by $\lambda$.  The composition of this linear functional on $V /
W$ with the canonical quotient mapping from $V$ onto $W$ gives a
linear functional on $V$ which extends $\lambda$ from $Z$ to $V$ and
has the same norm as $\lambda$ has on $Z$.

\section{Dual cones}
\label{dual cones}
\setcounter{equation}{0}

	Let $V$ be a finite-dimensional real vector space, and let
$\mathcal{C}$ be a closed convex cone in $V$.  To be more precise, one
can use an isomorphism between $V$ and ${\bf R}^n$ to define the
topology on $V$, i.e., so that the vector space isomorphism is a
homeomorphism.  This topology does not depend on the choice of
isomorphism with ${\bf R}^n$, because every invertible linear mapping
on ${\bf R}^n$ defines a homeomorphism from ${\bf R}^n$ onto itself.
Thus $\mathcal{C}$ is a closed subset of $V$ which contains $0$ and
has the property that $s \, v + t \, w \in V$ whenever $s, t$ are
nonnegative real numbers and $v, w \in \mathcal{C}$.

	Let us define $\mathcal{C}^*$ to be the set of linear
functionals $\lambda$ on $V$ such that $\lambda(v) \ge 0$ for all $v
\in \mathcal{C}$.  One can check that $\mathcal{C}^*$ defines a closed
convex cone in $V^*$.  This is called the dual cone associated to
$\mathcal{C}$.  It follows from the result in Section \ref{separation
theorems} for closed convex cones that $\mathcal{C}$ is actually equal
to the set of $v \in V$ such that $\lambda(v) \ge 0$ for all $\lambda
\in \mathcal{C}^*$.

	As a basic example, let $V$ be ${\bf R}^n$ for some positive
integer $n$, and let $\mathcal{C}$ be the closed convex cone
consisting of vectors $v = (v_1, \ldots, v_n)$ such that $v_j \ge 0$
for all $j$.  For each $w \in {\bf R}^n$, we get a linear functional
$\lambda_w$ on ${\bf R}^n$ by putting $\lambda_w(v) = \sum_{j=1}^n v_j
\, w_j$, and every linear functional on ${\bf R}^n$ arises in this
manner.  For this cone $\mathcal{C}$, the dual cone $\mathcal{C}^*$
consists of the linear functionals $\lambda_w$ such that $w \in
\mathcal{C}$, as one can readily verify.

	Now let $W$ be a finite-dimensional real or complex vector
space equipped with an inner product $\langle w, z \rangle$, and let
$V$ be the real vector space of linear mappings $A$ from $W$ to itself
which are self-adjoint, which is to say that $\langle A(w), z \rangle$
is equal to $\langle w, A(z) \rangle$ for all $w, z \in W$.  A
self-adjoint linear transformation $A$ on $W$ is said to be
nonnegative if $\langle A(w), w \rangle$ is nonnegative real number
for all $w \in W$, and the nonnegative self-adjoint linear
transformations on $V$ form a closed convex cone in the real vector
space of self-adjoint linear transformations on $V$.  If $T$ is a
self-adjoint linear transformation on $V$, then we get a linear
functional $\lambda_T$ on the vector space of self-adjoint linear
transformations on $W$ by setting $\lambda_T(A)$ equal to the trace of
$A \circ T$ for any self-adjoint linear transformation $A$ on $W$, and
every linear functional on the vector space of self-adjoint linear
transformations on $W$ arises in this manner.  If $T$ is a
self-adjoint linear transformation on $W$, then $\lambda_T(A) \ge 0$
for all nonnegative self-adjoint linear transformations $A$ on $W$ if
and only if $T$ is nonnegative.  This can be verified using the fact
that a self-adjoint linear transformation on $W$ can be diagonalized
in an orthonormal basis.

\section{Operator norms}
\label{operator norms}
\setcounter{equation}{0}

	Let $V$, $W$ be finite-dimensional vector spaces, both real or
both complex, and let $\mathcal{L}(V, W)$ be the vector space of
linear mappings from $V$ to $W$.  More precisely, this is a real
vector space if $V$, $W$ are real vector spaces and it is a complex
vector space if $V$, $W$ are complex vector spaces.  Notice that the
dual $V^*$ of $V$ is the same as $\mathcal{L}(V, {\bf R})$ when $V$ is
a real vector space and it is the same as $\mathcal{L}(V, {\bf C})$
when $V$ is a complex vector space.

	Suppose that $V$, $W$ are equipped with norms $\|v\|_V$,
$\|w\|_W$.  If $T$ is a linear mapping from $V$ to $W$, then there is
a nonnegative real number $k$ such that
\begin{equation}
	\|T(v)\|_W \le k \, \|v\|_V
\end{equation}
for all $v \in V$.  This can be derived from the case of mappings
between Euclidean spaces in the usual manner, and it basically amounts
to saying that $T$ is continuous as a mapping from $V$ to $W$ with
respect to the metrics associated to the norms on $V$, $W$.

	We define the operator norm of $T$ as a linear mapping
from $V$ to $W$ by
\begin{equation}
	\|T\|_{op, VW} = \sup \{\|T(v)\|_W : v \in V, \|v\|_V \le 1\}.
\end{equation}
Equivalently, $\|T\|_{op, VW}$ is the smallest nonnegative real number
that one can use as $k$ in the preceding paragraph.  It is easy to see
that $\|T\|_{op, VW}$ does indeed define a norm on $\mathcal{L}(V,
W)$.  In the case where $W$ is equal to ${\bf R}$ or ${\bf C}$, so
that $\mathcal{L}(V, W)$ is equal to $V^*$, the operator norm reduces
to the dual norm of a linear functional as defined in Section
\ref{dual spaces}.

	As another equivalent definition, the operator norm of $T$ is
equal to the supremum of $|\mu(T(v))|$ over $v \in V$, $\mu \in W^*$
where the $V$-norm of $v$ and the $W^*$-norm of $\mu$ are each less
than or equal to $1$.  That $|\mu(T(v))|$ is less than or equal to the
operator norm of $T$ for these $v$'s and $\mu$'s follows easily from
the definition.  The operator norm is equal to the supremum of these
quantities because norms in $W$ are detected by linear functionals as
in Section \ref{dual spaces}.

	For each linear mapping $T$ from $V$ to $W$ there is an
associated dual linear mapping $T^*$ from $W^*$ to $V^*$, defined by
saying that if $\mu$ is a linear functional on $W$, then $T^*(\mu)$ is
the linear functional on $V$ given by the composition of $\mu$ with
$T$.  One can check that the operator norm o $T^*$, with respect to
the dual norms on $V^*$, $W^*$, is equal to the operator norm of $T$
with respect to the original norms on $V$, $W$.

	As a special case, suppose that $V$ is equal to ${\bf R}^m$ or
${\bf C}^m$ for some positive integer $m$, equipped with the norm
$\|v\|_1 = |v_1| + \cdots + |v_m|$.  If $e_1, \ldots, e_m$ are the
standard basis vectors for $V$, so that the $l$th component of $e_j$
is equal to $1$ when $j = l$ and to $0$ when $j \ne l$, then one can
show that the operator norm of a linear mapping $T$ from $V$ into a
normed vector space $W$ is equal to the maximum of $\|T(e_1)\|_W,
\ldots, \|T(e_m)\|_W$.

	Suppose instead that $W$ is equal to ${\bf R}^n$ or ${\bf
C}^n$ for some positive integer $n$, equipped with the norm
$\|w\|_\infty = \max(|w_1|, \ldots, |w_n|)$.  A linear mapping $T$
from a normed vector space $V$ into $W$ is basically the same as a
collection $\lambda_1, \ldots, \lambda_n$ of $n$ linear functionals on
$V$, corresponding to the $n$ components of $T(v)$ in $W$.  The
operator norm of $T$ is then equal to the maximum of the dual norms of
$\lambda_1, \ldots, \lambda_n$ with respect to the given norm on $V$.

\section{Trace norms}
\label{trace norms}
\setcounter{equation}{0}

	Let $V$, $W$ be finite-dimensional vector spaces, both real or
both complex, and let $\|\cdot \|_V$, $\|\cdot \|_W$ be norms on $V$,
$W$, respectively.  Also let $T$ be a linear mapping from $V$ to $W$.
We can express $T$ as
\begin{equation}
	T(v) = \sum_{j=1}^l \lambda_j(v) \, w_j,
\end{equation}
for some linear functionals $\lambda_1, \ldots, \lambda_l$ on $V$ and
some vectors $w_1, \ldots, w_l$ in $W$.

	If $\lambda \in V^*$ and $w \in W$, then $\lambda(v) \, w$
defines a linear mapping from $V$ to $W$ with rank $1$, unless
$\lambda$ and $w$ are both $0$ in which event the linear mapping is
$0$.  The operator norm of this linear mapping is equal to the product
of the dual norm of $\lambda$ and the norm of $w$ in $W$.  If $T$ is
given as a sum as above, then for each $j$ one can take the product of
the dual norm of $\lambda_j$ and the norm of $w_j$ in $W$, and the sum
of these products is a nonnegative real number which is greater than
or equal to the operator norm of $T$.  The \emph{trace norm} of $T$ is
denoted $\|T\|_{tr, VW}$ and defined to be the infimum of this sum of
products over all such representations of $T$.  It is easy to see that
the trace norm does indeed define a norm.  In particular we have that
\begin{equation}
	\|T\|_{op, VW} \le \|T\|_{tr, VW}
\end{equation}
by the earlier remarks, which shows that the trace norm of $T$ is
equal to $0$ if and only if $T$ is equal to $0$.  At any rate, the
homogeneity and subadditivity of the trace norm can be derived
directly from the definition.

	If $T$ is a linear mapping from $V$ to $W$ and $A$ is a linear
mapping from $W$ to $V$, then the composition $A \circ T$ is a linear
mapping from $V$ to itself, and we can take its trace $\tr A \circ T$
in the usual manner.  If $T(v) = \lambda(v) \, w$ for some $\lambda
\in V^*$ and $w \in W$, then $(A \circ T)(v) = \lambda(v) \, A(w)$,
and the trace of $A \circ T$ is equal to $\lambda(A(w))$.  Using this
one can check that
\begin{equation}
	|\tr A \circ T| \le \|A\|_{op, WV} \, \|T\|_{tr, VW}
\end{equation}
for all linear mappings $A : W \to V$ and $T : V \to W$.  One can
think of $T \mapsto \tr A \circ T$ as a linear functional on
$\mathcal{L}(V, W)$, and the dual norm of this linear functional with
respect to the trace norm on $\mathcal{L}(V, W)$ is equal to
$\|A\|_{op, WV}$.

\section{Vector-valued functions}
\label{vector-valued functions}
\setcounter{equation}{0}

	Let $E$ be a finite nonempty set.  If $V$ is a real or complex
vector space, let us write $\mathcal{F}(E, V)$ for the vector space of
functions on $E$ with values in $V$.  Thus $\mathcal{F}(E, {\bf R})$,
$\mathcal{F}(E, {\bf C})$ denote the vector spaces of real or
complex-valued functions on $E$.

	If $1 \le p \le \infty$, let us write $\|f\|_p$ for the usual
$p$-norm of a real or complex-valued function $f$ on $E$, so that
\begin{equation}
	\|f\|_p = \bigg(\sum_{x \in E} |f(x)|^p \bigg)^{1/p}
\end{equation}
when $1 \le p < \infty$ and
\begin{equation}
	\|f\|_\infty = \max \{|f(x)| : x \in E\}.
\end{equation}
Let $V$ be a vector space equipped with a norm $\|v\|$, and let us
extend the $p$-norms to $V$-values functions on $E$ by putting
\begin{equation}
	\|f\|_{p, V} = \bigg(\sum_{x \in E} \|f(x)\|^p \bigg)^{1/p}
\end{equation}
when $1 \le p < \infty$ and
\begin{equation}
	\|f\|_{\infty, V} = \max \{\|f(x)\| : x \in E\},
\end{equation}
for any $V$-valued function $f$ on $E$.  It is easy to check that
these do define norms on $\mathcal{F}(E, V)$, using the properties of
the norm $\|v\|$ on $V$ and the $p$-norms for scalar-valued functions.

	Let $T$ be a linear transformation on the vector space of real
or complex-valued functions on $E$.  If $V$ is a real or complex
vector space, as appropriate, then $T$ induces a natural linear
transformation on the vector space of $V$-Valued functions on $E$.
For instance, for each $x \in E$ and scalar-valued function $f$ on
$E$, $T(f)(x)$ is a linear combination of $f(y)$, $y \in E$, and one
can use the same coefficients to define $T(f)$ when $f$ is a
vector-valued function.

	In general the relationship between the operator norm of $T$
acting on scalar-valued functions and the operator norm of $T$ acting
on vector-valued functions can be complicated.  If $V$ happens to be
${\bf R}^n$ or ${\bf C}^n$ equipped with a $p$-norm, and if we use the
$p$-norm for scalar valued functions on $E$, then the operator norm of
$T$ on scalar-valued functions with respect to the $p$-norm will be
the same as the operator norm of $T$ acting on $V$-valued functions
with respect to the norm $\|f\|_{p, V}$ defined above.

	Fix $p$, $1 \le p < \infty$, and a positive integer $n$.  Let
${\bf S}^{n-1}$ denote the standard unit sphere in ${\bf R}^n$
equipped with the Euclidean norm, which is the compact set of vectors
with Euclidean norm equal to $1$.  Let $V$ be the vector space of
real-valued continuous functions on ${\bf S}^{n-1}$.

	We can define a $p$-norm on $V$ by taking the integral of the
$p$th power of the absolute value of a continuous real-valued function
on ${\bf S}^{n-1}$, and then taking the $(1/p)$th power of the result.
Here we use the standard element of integration on ${\bf S}^{n-1}$
which is invariant under rotations, and which one may wish to
normalize so that the total measure of the sphere is equal to $1$.

	If $T$ is a linear transformation acting on real-valued
functions on $E$, then we get a linear transformation acting on
$V$-valued functions as before.  The operator norm of $T$ acting on
real-valued functions and using the $p$-norm on them is equal to the
operator norm of the associated linear transformation acting on
$V$-valued functions using the norm $\|f\|_{p, V}$ based on the
$p$-norm for functions on ${\bf S}^{n-1}$ described in the preceding
paragraph.

	Let $L$ denote the linear subspace of $V$ consisting of
functions on ${\bf S}^{n-1}$ which are restrictions of linear
functions from ${\bf R}^n$, so that each function in $L$ is given by
the inner product of the point in the sphere with some fixed vector in
${\bf R}^n$.  In this way we can identify $L$ with ${\bf R}^n$, and
the restriction of the $p$-norm on $V$ to $L$ corresponds to a
constant multiple of the usual Euclidean norm on ${\bf R}^n$.  Using
this it follows that if $T$ is a linear operator acting on real-valued
functions on $E$, then the operator norm of $T$ with respect to the
$p$-norm on real-valued functions on $E$ is the same as the operator
norm of the corresponding linear transformation acting on ${\bf
R}^n$-valued functions, using the norm on ${\bf R}^n$-valued functions
obtained from the $p$-norm on real-valued functions and the Euclidean
norm on ${\bf R}^n$.

\end{document}